\input amstex
\input amsppt.sty
\magnification\magstep1
\input xepsf
\def\ni\noindent
\def\sbs{\subset}

\def\as{\operatorname{asdim}}

\def\diam{\operatorname{diam}}
\def\dim{\operatorname{dim}}
\def\Int{\operatorname{Int}}

\def\ord{\operatorname{ord}}

\def\R{\text{\bf R}}

\def\N{\text{\bf N}}

\def\sV{\Cal V}

\def\sU{\Cal U}

\def\p{\partial}

\hoffset= 0.0in \voffset= 0.0in \hsize=32pc \vsize=38pc
\baselineskip=24pt \NoBlackBoxes \topmatter
\author
A. Dranishnikov
\endauthor

\title
On asymptotic dimension of amalgamated products and right-angled
Coxeter groups
\endtitle
\abstract We prove the inequality
$$
\as A\ast_CB\le\max\{\as A,\as B,\as C+1\}
$$
Then we apply this inequality to show that
the asymptotic dimension of any right-angled Coxeter group does not
exceed the dimension of its Davis' complex.
\endabstract

\thanks The author was partially supported by NSF grants DMS-0604494
\endthanks

\address University of Florida, Department of Mathematics, P.O.~Box~118105,
358 Little Hall, Gainesville, FL 32611-8105, USA
\endaddress


\email  dranish\@math.ufl.edu
\endemail

\endtopmatter

\document
\head \S0 Introduction \endhead

Asymptotic dimension was introduced by Gromov as an invariant of
finitely generated groups [Gr]. It is defined for metric spaces
and applied to finitely generated groups with the word metric.
Since by the definition it is quasi-isometry invariant, it does
not depend on the choice of a finite generating set. It turns out
that the asymptotic dimension is a coarse invariant in sense of
Roe [Ro]. Since all proper left invariant metrics on any countable
group are coarsely equivalent ([DSm], [Sh]), the notion of
asymptotic dimension can be extended to all countable groups.

The interest to the asymptotic dimension was sparked by Gouliang
Yu's proof of the Novikov Higher Signature conjecture for manifolds
whose fundamental group has finite asymptotic dimension [Yu1].
Similar progress on related conjectures was done under assumption of
finite asymptotic dimension in the works [Ba],[CG],[Dr1],[DFW].

Finite asymptotic dimensionality is proven for many classes of
groups. The exact computation of asymptotic dimension of groups is a
more difficult task. To the best of my knowledge it is completed
only for polycyclic groups and for hyperbolic groups. For polycyclic
groups the asymptotic dimension equals the Hirsch length:
$\as\Gamma=h(\Gamma)$ (see [BD3] for the inequality in one direction
and [DSm] for the other direction). The asymptotic dimension of a
finitely generated hyperbolic group equals the covering dimension of
its boundary plus one: $\as\Gamma=\dim\p_{\infty}\Gamma+1$ [B],[BL].
In view of Bestvina-Mess' formula [BM] we have $\as\Gamma=
vcd(\Gamma)$ for those hyperbolic groups for which the virtual
cohomological dimension is defined (e.g. for residually finite
hyperbolic groups).

The Coxeter groups are considered as a playground for many problems
and conjectures in Geometric Group Theory.  All Coxeter groups have
finite asymptotic dimension since they isometrically embeddable in
the finite product of trees [DJ]. The original such embedding is due
to Januszkiewicz and it gives an estimate $\as\Gamma\le|S|$ for a
Coxeter system $(\Gamma,S)$. In [DSc] we noticed that Januszkiewicz
technique can be pushed to bring a better estimate to $\as\Gamma\le
ch(N(\Gamma))$ where $N(\Gamma)$ is the nerve of $(\Gamma,S)$ and
$ch(N(\Gamma))$ is the chromatic number of the 1-skeleton of
$N(\Gamma)$. A low bound for the asymptotic dimension of Coxeter
groups is given by Corollary 4.11 in [Dr2]: $vcd\Gamma\le\as\Gamma$.
Since $\dim N+1\le ch(N)$ for every simplicial complex $N$ and
$vcd(\Gamma)\le\dim\Sigma(\Gamma)=\dim N(\Gamma)+1$, it is natural
to assume that $\as\Gamma\le\dim N(\Gamma)+1$ where $\Sigma(\Gamma)$
denotes the Davis complex of $(\Gamma,S)$. In this paper we prove
this inequality for right-angled Coxeter groups.  Generally, the
estimate $\as\Gamma\le\dim N(\Gamma)+1$ is not optimal. Perhaps the
most natural guess would be that $\as\Gamma= vcd\Gamma$ but in view
of Bestvina's candidate among right-angled Coxeter groups for a
counterexample to the Eilenberg-Ganea problem [D] this conjecture
seems to be very difficult. The proposed equality has to be checked
first for $\Gamma$ with 2-dimensional acyclic nerves $N(\Gamma)$.

In this paper the estimate $\as\Gamma\le\dim N(\Gamma)+1$ is
proven by induction on dimension of the nerve. The main ingredient
here is the inequality for the asymptotic dimension of the
amalgamated product
$$
(*)\ \ \ \ \ \ \ \as A\ast_CB\le\max\{\as A,\as B,\as C+1\}
$$
which is proven in this paper.

This inequality was conjectured in [BD2] as the equality. It was
proven when $C$ is a finite group in [BD3]. Later  J-P. Caprice
found an example among Coxeter groups where the inequality is strict.

As a corollary of (*) this paper gives a new proof of the free
product equality
$$ \as A\ast B=\max\{\as A,\as B,1\}.$$
The existed proof in [BDK] is quite long and it appeals to the
asymptotic inductive dimension theory developed in [DZ].

\head \S1 Asymptotic dimension \endhead

We recall the definition of asymptotic dimension of a metric space [Gr]:
{\it $\as X\le n$ if for every $r<\infty$ there exist uniformly bounded,
$r$-disjoint families $\sU^0,\ldots, \sU^n$ of subsets of $X$
such that $\cup_i\sU^i$ is a cover of $X$.}

Let $r\in\R_+$ be given and let $X$ be a metric space. We will say
that a family $\sU$ of subsets of $X$ is $r$-{\it disjoint} if
$d(U,U')\ge r$ for every $U\neq U'$ in $\sU.$ Here,
$d(U,U')=\inf\{d(x,x')\mid x\in U, x'\in U'\}.$

For a cover $\sU$  of a metric space $X$ we denote by
$L(\sU)=\inf_{U\in\sU}\sup_{x\in X}d(x,X\setminus U)$ the Lebesgue
number of $\sU$. We recall that the order $\ord\sU$ of a cover $\sU$
is the maximal number of elements with the nonempty intersection.

We say that $(r,d)$-$\dim X\le n$ if for every $r>0$ there exists a
$d$-bounded cover $\sU$ of $X$ with $\ord\sU\le n+1$ and with the
Lebesgue number $L(\sU)>r$. We refer to such a cover as to an
$(r,d)$-cover of $X$.

\proclaim{Proposition 1.1}[BD2] For a metric space $\as X\le n$ if
and only if there is a function $d(r)$ such that $(r,d(r))$-$\dim
X\le n$ for all $r>0$.
\endproclaim
Let $B_R(x)$ denote the closed $R$-ball centered at $x$ and let
$N_R(A)=\{x\in X\mid d(x,A)\le R\}$ denotes the closed
$R$-neighborhood of $A$. Thus, $B_R(x)=N_R(\{x\})$.

\proclaim{Proposition 1.2} Suppose that $X\subset Y$ is given the
restriction metric and let $\sU$ be an $(r,d)$-cover of $X$. Then
$N_{r/4}(X)$ admits an $(r/4,d+r)$-cover $\tilde{\sU}$ with
$\ord\tilde{\sU}\le ord\sU$.
\endproclaim
\demo{Proof} For every $U\in\sU$ we define $$\bar U=\bigcup\{Int
B_{r/2}(x)\mid d(x,X\setminus U)\ge r\}.$$ Clearly,
$N_{r/2}(X)\subset\cup_{U\in\sU}\bar U$. We show that $\ord\{\bar
U\mid U\in\sU\}\le \ord\sU$. Let $y\in \bar U_1\cap\dots\cap\bar
U_k$. Let $x_i\in U_i$ be such that $d(x_i,y)<r/2$ and
$d(x_i,X\setminus U_i)\ge r$. Since $d(x_i,x_1)<r$ and
$d(x_i,X\setminus U_i)\ge r$, it follows that $x_1\in U_i$ for all
$i$. Thus, $U_1\cap\dots\cap U_k\ne\emptyset$.

Let $\tilde U=\bar U\cap N_{r/4}(X)$. Then $\tilde{\sU}=\{\tilde
U\}$ is an $(r/4,d+r)$-cover $\tilde{\sU}$ with $\ord\tilde{\sU}\le
\ord\sU$. \qed
\enddemo

Let $K$ be a countable simplicial complex.  There is a metric on
$|K|$ called {\it uniform} which comes from the geometric realization
of $K$ .  It is defined by embedding of $K$ into the Hilbert space
$\ell^2=\ell^2(K^{(0)})$ by mapping each vertex $v\in K^{(0)}$ to
a corresponding element of an orthonormal basis for $\ell^2$ and giving $K$ the
metric it inherits as a subspace.

A map $\varphi:X\to Y$ between
metric spaces is {\it uniformly cobounded} if for every $R>0$,
$\diam(\varphi^{-1}(B_R(y)))$
is uniformly bounded.
We call a map $\varphi:X\to |K|$ to a simplicial complex {\it $c$-cobounded},
$c\in \R_+$, if $\diam(\varphi^{-1}(\Delta))< c$ for all
simplices $\Delta\subset K$.

The following was proven by Gromov [Gr] (see also [BD2], [Ro]).
\proclaim{Theorem 1.1} Let $X$ be a metric space.  The following
conditions are equivalent. \roster
    \item $\as X\le n$;
    \item for every $\epsilon>0$ there is a uniformly cobounded,
    $\epsilon$-Lipschitz map $\varphi:X\to K$ to a uniform simplicial
    complex of dimension $n.$
\endroster
\endproclaim

This theorem is proved by using projections to the nerves of open covers.
The projection $p_{\sU}:X\to Nerve(\sU)\subset\ell^2(\sU)$ defined by the formula
$$
p_{\sU}(x)=(\phi_U)_{U\in\sU}, \ \ \ \ \ \ \
\phi_U(x)\frac{d(x,X\setminus U)}{\sum_{V\in \sU} d(x,X\setminus V)}
$$
is called {\it canonical}.

A map $f:X\to
Y$ between metric spaces is a {\it coarse embedding} if there exist
non-decreasing functions $\rho_1$ and $\rho_2$,
$\rho_i:{\R}_+\to{\R}_+$ such that
$\rho_i\to\infty$ and for every $x,x'\in X$
$$\rho_1(d_X(x,x'))\le d_Y(f(x),f(x'))\le
\rho_2(d_X(x,x'))\text{.}$$ Such a map is often called a {\it
coarsely uniform embedding} or just a {\it uniform embedding}. The
metric spaces $X$ and $Y$ are {\it coarsely equivalent} if there is
a coarse embedding $f:X\to Y$ so that there is some $R$ such that
$Y\subset N_R(f(X))$.

Observe that quasi-isometric spaces are coarsely equivalent with
linear $\rho_i.$

\proclaim{Proposition 1.3} Let $f:X\to Y$ be a coarse equivalence.
Then $\as X=\as Y.$
\endproclaim
As a corollary we obtain that $\as\Gamma$ is an invariant for
finitely generated groups. One can extend this definition of $\as$
for all countable groups by considering left-invariant proper
metrics on $\Gamma$. All such metrics are coarsely equivalent
[DSm], [Sh].

\proclaim{Theorem 1.2} {\rm [DSm]} Let $G$ be a countable group.
Then $\as G=\sup\as F$ where the supremum is taken over all
finitely generated subgroups $F\subset G$.
\endproclaim

For a subset $Y\subset X$ of metric space $X$ when we write $\as Y$ we
assume that $Y$ is taken with the metric obtained by restriction.

Also we use in this paper the following two theorems from [BD1]:
\proclaim{Finite Union Theorem} For every metric space presented as
a finite union $X=\cup X_i$ there is the formula
$$\as(\cup X_i)=\max\{\as X_i\}.$$
\endproclaim

\proclaim{Infinite Union Theorem} Let $X=\cup_\alpha X_\alpha$ be
a metric space where the family $\{X_\alpha\}$ satisfies the
inequality $\as X_\alpha\le n$ uniformly.  Suppose further that
for every $r$ there is a $Y_r\subset X$ with $\as Y_r\le n$ so
that $d(X_\alpha\setminus Y_r,X_{\alpha'}\setminus Y_r)\ge r$
whenever $X_\alpha\neq X_{\alpha'}.$  Then $\as X\le n.$
\endproclaim
We recall that the family
$\{X_\alpha\}$ of subsets of $X$ satisfies the inequality $\as
X_\alpha\le n$ {\it uniformly} if for every $r<\infty$ one can find
a constant $R$ so that for every $\alpha$ there exist $r$-disjoint
families $\sU^0_\alpha,\ldots,\sU^n_\alpha$ of $R$-bounded subsets
of $X_\alpha$ covering $X_\alpha.$

The following Propositions are taken from [BD2] (Proposition 2 and
Lemma 1).

\proclaim{Proposition 1.4} For every simplicial map $g:X\to Y$ the
mapping cylinder $M_{g}$ admits a triangulation with the set of
vertices equal to the disjoint union of vertices of $X$ and $Y$.
\endproclaim

We consider the uniform metric on $M_g$.

For a cover $\sU$  of a metric space $X$ we denote by
$b(\sU)=\sup_{U\in\sU}diam(U)$ the diameter of $\sU$. We note that
if for two covers $b(\sV)< L(\sU)$ then there is a map $G:\sV\to\sU$
with the property $G(V)\subset U$. Note that any such map
$G:\sV\to\sU$ defines a simplicial map $g:Nerve(\sV)\to Nerve(\sU)$
of the nerves.

We use the notations $\p N_r(A)=\{x\mid d(x,A)=r\}$ for the boundary
of the $r$-neighborhood and $r$-$\Int(A)=A\setminus N_r(X\setminus
A)$ for the {\it $r$-interior} of $A$.

\proclaim{Lemma 1.1} [BD2] For every $n\in N$ there is a monotone
tending to infinity function $\mu:\R_+\to\R_+$ with the following
property: Given $\epsilon>0$, let $W\subset X$ be a subset of a
geodesic metric space $X$ and let $\lambda\ge 1/\epsilon$. Then for
every two covers $\sV$ of $N_{\lambda}(\p W)$ and $\sU$ of $W$
by open subsets of $X$ with
the order $\le n+1,$ and with $L(\sU)>b(\sV)>L(\sV)\ge\mu(\lambda)$,
there is a $2b(\sU)$-cobounded $\epsilon$-Lipschitz map $f: W\to
M_g$ to the mapping cylinder of a simplicial map $g:Nerve(\sV)\to
Nerve(\sU)$ between the nerves  such that $f|_{\p W}=p_{\sV}|_{\p
W}$ where $p_{\sV}:N_{\lambda}(\p W)\to Nerve(\sV)\subset M_g$ is
the canonical projection.
\endproclaim

We note that the formulation of Lemma 1.1 differs slightly from
Lemma 1 in [BD2]. Namely, Lemma 1.1 becomes Lemma 1 if one consider
the case $W=N_r(\p W)$. Nevertheless the same formula for $f$ and
the same proof as in [BD2] are valid for the general case.

A {\it partition} of a metric space $X$ is a presentation as a union
$X=\cup_iW_i$ such that $\Int(W_i)\cap \Int(W_j)=\emptyset$ whenever
$i\ne j$.

\proclaim{Partition Theorem} Let $X$ be a geodesic metric space.
Suppose that for every $R>0$ there is $d>0$ and a partition
$X=\cup_i^{\infty}W_i$ with  $\as W_i\le n$ uniformly on $i$ such
that $(R,d)$-$\dim (\cup_i\partial W_i)\le n-1$ where
$\cup_i\partial W_i$ is taken with a metric restricted from $X$.
Then $\as X\le n$.
\endproclaim
\demo{Proof} We apply Theorem 1.1. Given $\epsilon>0$ we construct a
uniformly cobounded $\epsilon$-Lipschitz map $\phi:X\to K$. We apply
the assumption with $R=4\mu(1/\epsilon)$ where $\mu$ is taken from
Lemma 1.1. Let $\lambda=1/\epsilon$. Since $\lambda\le R/4$ and
$\mu(t)\ge t$, by Proposition 1.2 there is an $(r/4,2d)$-cover $\sV$
of $N_{\lambda}(\cup_i\partial W_i)$ of order $\le n$. We may assume
that it is a cover by open in $X$ sets. Let $\sV_i=\sV|_{\p W_i}$ be
the restriction, i.e., $\sV_i$ consists of those elements of $\sV$
that have a nonempty intersection with $\p W_i$. Let $\sU_i$ be a
cover of $W_i$ with $L(\sU_i)>2d\ge b(\sV_i)$ and with $b(\sU_i)<D$
for all $i$ for some fixed $D$. By Lemma 1.1 there is a
$2D$-cobounded $\epsilon$-Lipschitz map $f_i:W_i\to M_{g_i}$ to a
uniform complex where $M_{g_i}$ is the mapping cylinder of a
simplicial map $g_i:Nerve (\sV_i)\to Nerve(\sU_i)$ and $f_i$
coincides on $\p W_i$ with the canonical projection to the nerve
$p_{\sV}:\cup_i\partial W_i\to Nerve(\sV)$. We define
$$K=(Nerve(\sV)\coprod\coprod_iM_{g_i})/\sim$$
as the quotient space under identification along the complexes
$Nerve(\sV_i)$. Then the union of $f_i$ defines a map $f:X\to K$.
Clearly, $f$ is $2D$-cobounded. Since $X$ is geodesic and each $f_i$
is $\epsilon$-Lipschitz, $f$ is $\epsilon$-Lipschitz with respect to
the uniform metric on $K$.\qed
\enddemo

\head \S2 Asymptotic dimension of amalgamated product
\endhead
Let $A$ and $B$ be finitely generated groups and let $C$ be a common
subgroup. We fix finite generating symmetric sets $S_A$ and $S_B$.
Let $d$ denote the word metric on $A\ast_CB$ corresponding to the
generating set $S_A\cup S_B$. The group $G=A\ast_CB$ acts on the
Bass-Serre tree whose vertices are the left cosets $G/A\coprod G/B$
and the vertices $xA$ and $xB$, $x\in G$, and only them are joined
by edges. The edges $[xA,xB]$ are labeled by the cosets $xC$. We
consider the action of $G$ on the dual graph $K$. Thus vertices of
$K$ are the left cosets $xC$. Two vertices $xC$ and $x'C$ are joined
by an edge if an only if the edges in the Bass-Serre tree with these
labels have a common vertex. Note that $K$ is a tree-graded space in
the sense of Drutu-Sapir [DS] with pieces $\Delta(A)$ and
$\Delta(B)$, the 1-skeletons of the simplices spanned by $A/C$ or
$B/C$. Thus, $K$ is partitioned into these pieces in a way that
every two pieces have at most one common vertex and the nerve of the
partition is a tree. The graph $K$ has an additional property that
all vertices are the intersection points of exactly two pieces of
the different types $\Delta(A)$ and $\Delta(B)$. We consider the
simplicial metric on $K$, i.e., every edge has length one and we use
the notation $|u,v|$ for the distance between vertices $u,v\in
K^{(0)}$. For $u\in K^{(0)}$ by $|u|$ we denote the distance to the
the vertex with label $C$. Note $K$ has the unique geodesic property
for every pair of vertices. There is a natural projection $\pi:G\to
K$ defined by the action: $\pi(g)=gC$.

\proclaim{Assertion 2.1} The map $\pi:G\to K$ extends to a
simplicial map of the Cayley graph of $G$, $\pi:C(G)\to K$.
\endproclaim
\demo{Proof} Let $g\in G$ and $s\in S_A\cup S_B$. If $s\in C$, then
$\pi(g)=\pi(gs)$ and the edge $[g,gs]\subset C(G)$ is mapped to the 
vertex $gC=gsC$. With out loss of generality we may assume that 
$s\in A\setminus C$. We need to show that
$\pi(g)=gC$ and $\pi(gs)=gsC$ are joined by an edge in $K$.  Note that $gA$
is the common vertex for the edges $[gB,gA]$ and $[gA,gsB]$ in the
Bass-Serre tree. Hence the vertices corresponding these edges are
joined by an edge in $K$. Thus, the vertices $gC$ and $gsC$ are
joined by an edge in $K$. \qed
\enddemo
As a corollary we obtain that $\pi$ is 1-Lipschitz.
\medskip
\epsfysize=3in \centerline{\epsffile{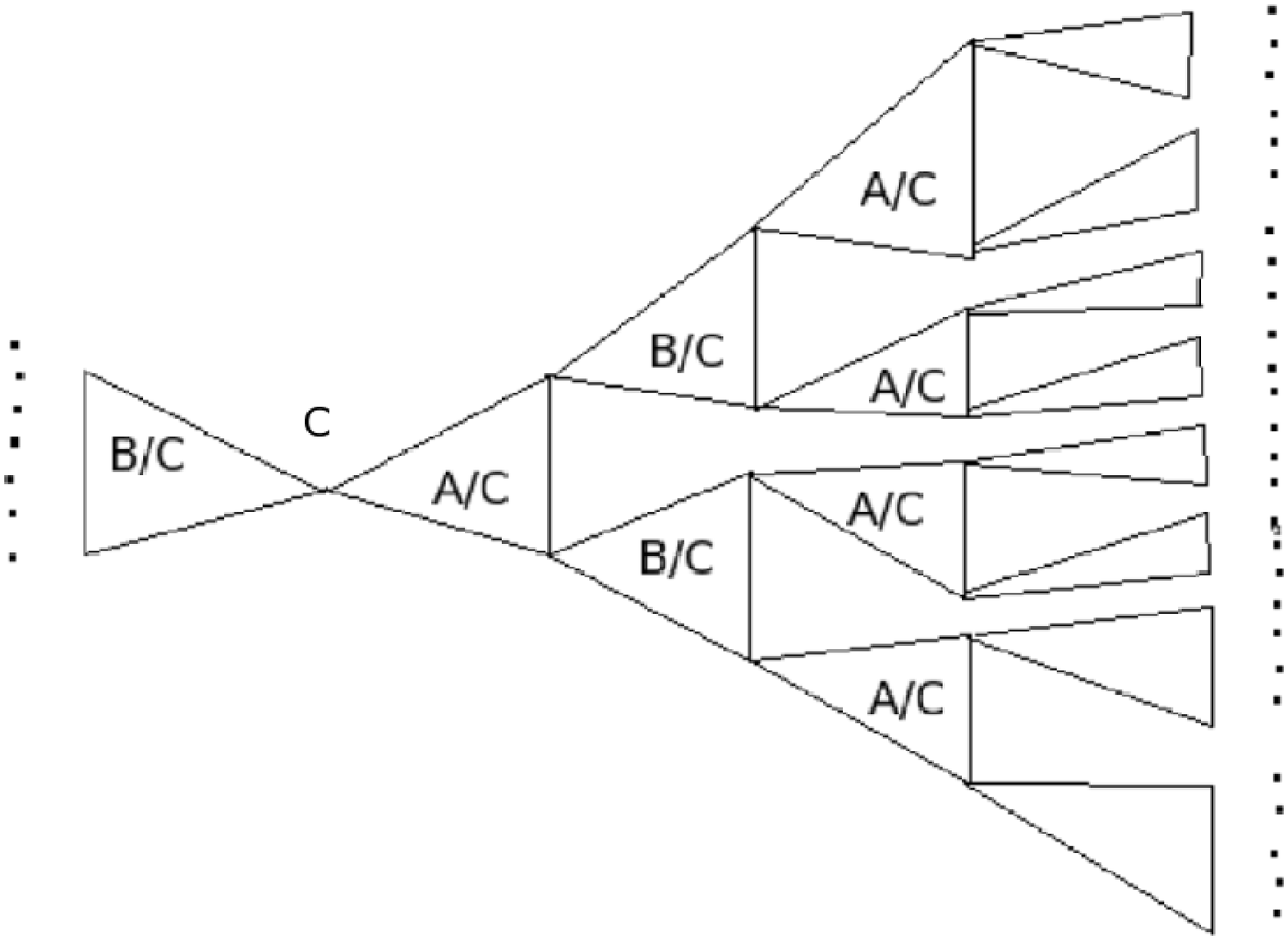}}
\smallskip
\centerline{The tree-graded complex $K$}
\medskip

The base vertex $C$ separates $K$ into two parts $K_A\setminus\{C\}$ 
and $K_B\setminus\{C\}$. Let $\bar d$ denote the graph metric on $K$.
We denote by $B_r^A$ the $r$-ball, $r\in\N$, in  $K_A$ centered at $C$.
There is a partial order on vertices of $K$ defined as follows:
$v\le u$ if and only if $v$ lies in the geodesic segment $[C,u]$
joining the base vertex with $u$. For $u\in K^{(0)}$ of nonzero
level and $r>0$ we denote by
$$K^u=\{v\in K^{(0)}\mid v\ge u\},\ \ \ B_r^u=\{v\in K^u\mid  |v|\le|u|+r\}.$$
For every vertex $u\in K^{(0)}$ represented by a coset $gC$ with
$|u|$ even we have the equalities $B^u_r=gB_r^A$, $K^u=gK_A$ and
hence the equalities $\pi^{-1}(B^u_r)=g\pi^{-1}(B_r^A)$ and
$\pi^{-1}(K^u)= g\pi^{-1}(K_A)$.

We say that a set $F\subset G$ {\it separates} two subsets $H_1$ and
$H_2$ in $G$ if it separates them in the Cayley graph $C(G)$, that
is every path in  $C(G)$ with the endpoints in $H_1$ and $H_2$ meets
$F$.

Let $D_R=\{x\in G\mid d(x,C)=R\}$ be the boundary of the
$R$-neighborhood of $C$ in $\pi^{-1}(K_A)$, $R\in\N$. For $u\in
K^{(0)}$ we denote $D_R^u=g_u(D_R)$, $g_u\in G$, where the coset
$g_uC$ represents $u$. Note that $\pi(D^u_R)\subset B^u_R$.

\proclaim{Proposition 2.1} For every vertex $u\in K$  with even
$|u|$ and for every $v\in K^{(0)}$  incomparable with $u$ or
satisfying $v<u$, the set $D_R^u$ separates $\pi^{-1}(v)$ and
$\pi^{-1}(u')$ in $A\ast_CB$ whenever $u<u'$ and $|u'|-|u|>R$.
\endproclaim
\demo{Proof}  Since $G$ acts by isometries, it suffices to show that
$D_R$ separates $\pi^{-1}(K_B)$ and $\pi^{-1}(u')$ with
$|u'|>R$ and $u'\in K_A$. In view of the fact that $\pi$ is
1-Lipschitz, $d(C,\pi^{-1}(u'))>R$. Hence  $D_R$ separates $C$ and
$\pi^{-1}(u')$. By Assertion 2.1  the image of a 
path in $C(G)$ is a  path in $K$. Since
every path in $K$ from $K_A$ to $K_B$ hits the vertex $C$, it
follows that every path in the Cayley graph from $\pi^{-1}(K_A\setminus C)$ to
$\pi^{-1}(K_B\setminus C)$ hits the set $C=\pi^{-1}(\{C\})$. Hence every path from
$\pi^{-1}(u')$ to $\pi^{-1}(K_B)$ hits $D_R$.\qed
\enddemo

\proclaim{Proposition 2.2} If $R\le r/4$, then $d(gD_R, g'D_R)\ge
2R$ for $g,g'\in G$ with $|gC|,|g'C|\in nr$, $n\in\N$, and $gC\ne
g'C$.
\endproclaim
\demo{Proof}  
We use notations $u=gC$ and $u'=g'C$ for the vertices in $K$.
First consider the case when $|u|\ne |u'|$.
Since
$\pi(gD_R)\subset B^u_R$, $\bar d(B^u_R,B^{u'}_R)\ge r-R\ge 3R$, and
$\pi$ is 1-Lipschitz, we obtain that $d(gD_R,g'D_R)\ge 3R$.

Let $|u|=|u'|$ and let $x\in gD_R$ and $y\in g'D_R$. Since every
path in $K$ between $\pi(x)$ and $\pi(y)$ goes through the vertices
$u$ and $u'$ and a path in the Cayley graph $C(G)$ is projected to
a path in $K$ (see Assertion 2.1), a geodesic from $x$ to $y$ in $C(G)$
passes through $gC$ and $g'C$. Since $d(x,gC)=R$ and
$d(y,g'C)=R$ we obtain the inequality
$d(x,y)\ge 2R$.\qed
\enddemo

We fix two set-theoretic sections $s_A:A/C\to A$ and $s_B:B/C\to B$
of $\pi_A$ and $\pi_B$ and denote by $X=im(s_A)\setminus C$ and
$Y=im(s_B)\setminus C$. These sections give rise the normal
presentation of elements in $A\ast_CB$. Namely, every element
$\gamma\in A\ast_CB$ can be presented uniquely in the following form
$\gamma=z_1\dots z_kc$ where $c\in C$, $z_i\in X\cup Y$, and $z_i$
are alternating in a sense that if $z_i\in X$ then $z_{i+1}\in Y$
and if $z_i\in Y$ then $z_{i+1}\in X$. Denote by $l(z_1\dots
z_kc)=k$ the length of the normal presentation. Clearly,
$l(\gamma)=|\gamma C|$ where $\gamma C$ is treated as a vertex of
$K$. If $X'$ and $Y'$ is a different choice of representatives and
$\gamma=z_1\dots z_kc$ and $\gamma=z_1'\dots z_k'c'$ are
corresponding normal presentations, then $z_1'=z_1c_1$,
$z_2'=c_1^{-1}z_2c_2$, $\ \dots \ $, $z_i'=c_{i-1}^{-1}z_ic_i$, $\
\dots\ $, $z_{k-1}'=c_{k-2}^{-1}z_{k-1}c_{k-1}$, and
$z_k'c'=c_{k-1}^{-1}z_kc$ where $c_i\in C$, $i=1,\dots, k$.

\proclaim{Assertion 2.2} Let $\gamma\in A\ast_CB$. Then
$\|\gamma\|\ge d(\beta_kc,C)$ for the normal presentation
$\gamma=\beta_1\dots\beta_kc$ for any choice of representatives $X$
and $Y$.
\endproclaim
\demo{Proof}  Let $\gamma=t_1\dots t_n$, $n=\|\gamma\|$, be the
shortest presentation. Then the word $t_1\dots t_n$ can be
partitioned into a normal presentation (for some choice of $X$ and
$Y$) $\gamma=\alpha_1\dots\alpha_k$, $\alpha_i\in A$ (or
$\alpha_i\in B)$, with $k=l(\gamma)$. Then $\|\gamma\|\ge
\|\alpha_k\|= d(\alpha_k,1)\ge d(\alpha_k,C)$.  If $X'$ and $Y'$ is
a different choice of representatives and $\gamma=\beta_1\dots
\beta_kc$ are corresponding normal presentations, then
$\beta_1=\alpha_1c_1$, $\beta_2=c_1^{-1}\alpha_2c_2$, $\ \dots \ $,
$\beta_i=c_{i-1}^{-1}\alpha_ic_i$, $\ \dots\ $,
$\beta_{k-1}=c_{k-2}^{-1}\alpha_{k-1}c_{k-1}$, and
$\beta_kc'=c_{k-1}^{-1}\alpha_k$ where $c_i\in C$, $i=1,\dots, k$.

Thus, $\beta_kc'=c\alpha_k$.  This implies that $\|\gamma\|\ge
d(\beta_kc,C)$. \qed
\enddemo

\proclaim{Lemma 2.1} Let $\as A,\as B\le n$. Then $\as (AB)^m\le n$
for all $m$ where $(AB)^m=AB\dots AB\subset A\ast_CB$.
\endproclaim
\demo{Proof} We prove that $\as AB\dots A(B)\le n$ by induction on
the length of the product $k$. The inequality is a
true statement for $k=1$. Assume that it holds for $k$. For the
sake of concreteness assume that $k$ is odd.
Thus, $\as F_1\dots F_k\le n$ where $F_{2i-1}=A$ and $F_{2i}=B$.
We show that $\as F_1\dots
F_kB\le n$. Consider the family $\{wB\mid
l(w)=k\}$. Since all sets $wB$ are isometric to $B$, $\as
wB\le k$ uniformly.

Given $r$ we define $Y_r=AB\dots ACB_r$ where $B_r$ is the $r$-ball
in $B$. We show that $d(wB\setminus Y_r,w'B\setminus Y_r)\ge r$. Let
$b,b'\in B\setminus CB_r$. Then $d(b,C)\ge r$ and $d(b',C)\ge r$.
Since $w^{-1}w'\notin B$, the normal presentation of
$b^{-1}w^{-1}w'b'$ ends with $c'b'$ for $c'\in C$. Then by Assertion
2.2, $\| b^{-1}w^{-1}w'b'\|\ge d(c'b',C)\ge r$. Then by the Infinite
Union Theorem we obtain that $\as (F_1\dots F_k\cap L_k)B\le n$
where $L_k$ is the set of all elements $w\in A\ast_CB$ with
$l(w)=k$, $F_i=A(B)$ are alternating, and $F_0=A$. Let $L_{<k}$ be
the set of all elements $w\in A\ast_CB$ with $l(w)<k$. The
inequality $\as (F_1\dots F_m\cap L_{<m})B\le n$ follows from
induction assumption and the Finite Union Theorem. \qed
\enddemo

This Lemma first appeared in [BD1]. We present it here with a proof
since the argument in [BD1] contains a gap in the proof of
Proposition 3.

 \proclaim{Theorem 2.1} For any finitely
generated groups $A$ and $B$ and a common subgroup $C$ there is the
inequality
$$
\as A\ast_C B\le\max\{\as A,\as B,\as C+1\}.
$$
\endproclaim
\demo{Proof} Let $n=\max\{\as A,\as
B,\as C+1\}$ and let $\pi:A\ast_CB\to
K=K_A\cup K_B$ be the projection to the graph dual to the Bass-Serre tree.
In view of the Finite Union Theorem it suffices to show that
$\as\pi^{-1}(K_A)\le n$ and $\as\pi^{-1}(K_B)\le n$. We prove the
first.

We apply the Partition Theorem. Let $R>0$ be given. Take $r>4R$.

In view of Proposition 2.1 $G=X_+\cup X_-$ with $X_+\cap X_-= D_R$
such that $X_+\subset\pi^{-1}(K_A)$, $\pi^{-1}(K_B)\subset X_-$ and
$D_R$ separates $X_+\setminus D_R$ and $X_-\setminus D_R$. For every
vertex $u\in K_A$ we fix an element $g_u\in G$ such that the coset
$g_uC$ represents $u$. We denote by $X_{\pm}^u=g_u(X_{\pm})$ and
define $V_r=X_+\cap(\bigcap_{|u|=r}X^u_-)$. Note that
$\pi(V_r)\subset B_{r+R}$. Let $V_r^u$ denote $g_u(V_r)$.  
Consider the partition
$$
\pi^{-1}(K_A)=\bigcup_{|u|=nr, n\in\N_+}V^u_r\cup N_R^A(C)$$ where
$N^A_R(C)=N_R(C)\cap\pi^{-1}(K_A)$.

Clearly, if $V^u_r\cap V_r^w\ne\emptyset$, then either $u<w$ and
$|w|=|u|+r$ or $w<u$ and $|u|=|w|+r$. If $V^u_r\cap
V_r^w\ne\emptyset$ and $u<w$ then $V^u_r\cap V_r^w= D_R^w$ where
$D_R^w=g_wD_R$.

Thus
$$
Z=\bigcup_{|u|=nr, n\in \N_+}\p V^u_r=\bigcup_{|u|=nr, n\in
\N_+}D_R^u.$$

We show that $(R,d)$-$\as Z\le n-1$ for some $d>0$. Since $D_R$ is
coarsely equivalent to $C$, we have $\as D_R\le n-1$. Hence there is
$d>0$ and an $(R,d)$-cover $\sU$ of $D_R$ with $\ord\sU\le n$. In
view of Proposition 2.2, $\tilde{\sU}=\cup_{|u|=nr, n\in\N_+}
g_u(\sU)$ is an $(R,d)$-cover of $Z$.

Since $\pi^{-1}(B_s)\subset (AB)^{s+1}$, by Lemma 2.1  we have
$\as\pi^{-1}(B_s)\le n$. Hence $\as\pi^{-1}(B_{r+R})\le n$ and
therefore, $\as V_r^u\le n$ uniformly. Note that $\as N_R^A(C)\le n$
(in fact, it is $\le n-1$).

By the Partition Theorem, $\as \pi^{-1}(K_A)\le n$\qed
\enddemo

\

\head \S3 Asymptotic dimension of right-angled Coxeter groups
\endhead

\

A symmetric matrix $M=(m_{ss'})_{s,s'\in S}$ is called a Coxeter
matrix if $m_{ss'}\in\N\cup\{\infty\}$, and $m_{ss}=1$ for all
$s\in S$. A group $\Gamma$ with a generating set $S$ is called a
Coxeter group if there is a Coxeter matrix $M=(m_{ss'})_{s,s'\in
S}$ such that $\Gamma$ admits a presentation
$$
\langle S\ \mid (ss')^{m_{ss'}}, s,s'\in S\rangle.
$$
A Coxeter group $\Gamma$ is called {\it even} if all finite
non-diagonal entries of $M$ are even. It is called {\it
right-angled} if all finite non-diagonal entries equal 2.

Every subset $W\subset S$ defines a subgroup $\Gamma_W\subset
\Gamma$ which will be called {\it parabolic}. Let $(\Gamma,S)$ be a
Coxeter group with a generating set $S$ and with a presentation
given by means of a Coxeter $S\times S$ matrix $M$. {\it The nerve}
$N(\Gamma)$ is a simplicial complex with the set of vertices $S$
where a subset $W\subset S$ spans a simplex if and only if the group
$\Gamma_W$ is finite. Thus, $s,s'\in S$, $s\ne s'$, form an edge if
and only if  $m_{ss'}\ne\infty$. We call the number $m_{ss'}$ a {\it
label} of the edge $[ss']$. By $N'$ we denote the barycentric
subdivision of $N$. The cone $C=Cone N'$ over $N'$ is called a {\it
chamber} for $\Gamma$. {\it The Davis complex}
$\Sigma=\Sigma(\Gamma,S)$ is the image of a simplicial map
$q:\Gamma\times C\to \Sigma$ defined by the following equivalence
relation on the vertices: $a\times v_{\sigma}\sim b\times
v_{\sigma}$ provided $a^{-1}b\in\Gamma_{\sigma}$ where $\sigma$ is a
simplex in $N$ and $v_{\sigma}$ is the barycenter of $\sigma$. We
identify $C$ with the image $q(e\times C)$. The group $\Gamma$ acts
properly and simplicially on $\Sigma$ with the orbit space
equivalent to the chamber. Thus, the Davis complex is obtained by
gluing the chambers $\gamma C$, $\gamma\in\Gamma$ along their
boundaries. The main feature of $\Sigma$ is that it is contractible
(see [D]).

 \proclaim{Theorem 3.1} For every right-angled Coxeter group
$\as\Gamma\le \dim N(\Gamma)+1$.
\endproclaim
\demo{Proof} We prove this inequality by induction on the
dimension of the nerve $N=N(\Gamma)$.  If $\dim N=0$,
then $\Gamma$ is virtually free group (possibly of zero rank) and
hence $\as\Gamma\le 1=\dim N+1$.

Let $\dim N=n$ and let $N$ be finite. We prove the inequality
$\as\Gamma\le n+1$ by induction on
the number of vertices in $N$. If this number is minimal, i.e.,
$n+1$, the inequality holds since the group $\Gamma$ in this case
is finite. We assume
that there is a vertex $v\in N$ such that the star $st(v,N)$ does
not contain all other vertices of $N$. If there is no such $v$, then
the 1-dimensional skeleton $N^{(1)}$ coincides with the 1-skeleton
of a simplex $\Delta$. Since the group is right-angled, $N=\Delta$,
and the group is finite,
so this case has been already considered. We take $K$ to be the link
$Lk(v,N)$ of such vertex $v$ and take $N_1$ to be the star $st(v,N)$
of this vertex. We define $N_2=N\setminus Ost(v,N)$ where $Ost(v,N)$ is
the open star of $v$.

Then $\Gamma=\Gamma_{N_1}\ast_{\Gamma_K}\Gamma_{N_2}$. By induction
assumption $\as \Gamma_K\le n$.
By the internal induction $\as\Gamma_{N_i}\le n+1$, $i=1,2$.
Then Theorem 2.1 implies that $\as\Gamma\le n+1$.
\qed
\enddemo
\proclaim{Corollary 3.1} For every right-angled Coxeter group
$\as\Gamma\le \dim \Sigma(\Gamma)$ where $\Sigma(\Gamma)$ is the
Davis complex.
\endproclaim
In view of recent result of Dymara and Schick [DySc] we obtain

\proclaim{Corollary 3.2} For a right-angled building $X$,
$\as X\le \dim X$.
\endproclaim

REMARK. In order to extend Theorem 3.1 to all Coxeter groups
one needs to show the inequality $\as\Gamma\le\dim N+1$
in the case when $N^{(1)}$ is the 1-skeleton of a simplex.

\enddocument